\documentclass[10pt,twoside]{article}
\usepackage{graphicx}
\usepackage{amsmath}
\usepackage{latexsym}
\usepackage{amsfonts}
\usepackage{amsthm}
\usepackage{Latex-document}

\newtheorem{theorem}{Theorem}[section]
\newtheorem{lemma}[theorem]{Lemma}
\newtheorem{conjecture}[theorem]{Conjecture}

\theoremstyle{definition}

\newtheorem{remark}[theorem]{Remark}

\def\doct{\delta_{oct}}
\def\dtet{\delta_{tet}}
\def\pt{\hbox{\it pt}}
\def\Vol{\hbox{vol}}
\def\R{{\mathbb R}}

\title{\bf  A Computer Verification of \vskip -2mm
the Kepler Conjecture\vskip 6mm}

\renewcommand{\thesection}{\arabic{section}}
\renewcommand{\thesubsection}{\thesection.\arabic{subsection}}

\def\addsec{\addtocounter{section}{1}}

\author{Thomas C. Hales\vspace*{-0.5cm}\thanks{Department of Mathematics,
University of Pittsburgh, Thackeray Hall, Pittsburgh PA 15260,
USA. E-mail: hales@pitt.edu}}
\date{\vspace{-8mm}}

\begin{document}

\maketitle

\thispagestyle{first} \setcounter{page}{795}

\begin{abstract}

\vskip 3mm

The Kepler conjecture asserts that the density of a packing of
congruent balls in three dimensions is never greater than
$\pi/\sqrt{18}$.  A computer assisted verification confirmed this
conjecture in 1998.   This article gives a historical introduction
to the problem. It describes the procedure that converts this
problem into an optimization problem in a finite number of
variables and the strategies used to solve this optimization
problem.

\vskip 4.5mm


\end{abstract}

\vskip 12mm

\section*{1. Historical introduction} \label{section 1}\setzero \addsec
\vskip-5mm \hspace{5mm }

\begin{figure}[htb]
  \centering
  \includegraphics{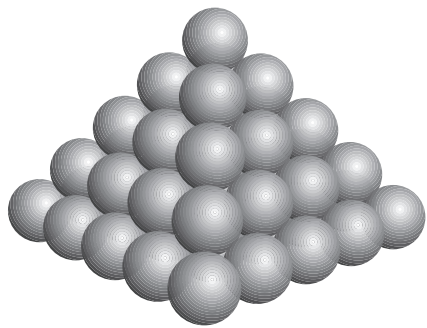}
  \caption{The face-centered cubic packing}
  \label{fig:fcc}
\end{figure}

The Kepler conjecture asserts that the density of a packing of
congruent balls in three dimensions is never greater than
$\pi/\sqrt{18}\approx 0.74048\ldots$. This is the oldest problem
in discrete geometry and is an important part of Hilbert's 18th
problem. An example of a packing achieving this density is the
face-centered cubic packing (Figure \ref{fig:fcc}).

A packing of balls is an arrangement of nonoverlapping balls of
radius 1 in Euclidean space. Each ball is determined by its
center, so equivalently it is a collection of points in Euclidean
space separated by distances of at least 2. The density of a
packing is defined as the $\limsup$ of the densities of the
partial packings formed by the balls inside a ball with fixed
center of radius $R$. (By taking the $\limsup$, rather than
$\liminf$ as the density, we prove the Kepler conjecture in the
strongest possible sense.) Defined as a limit, the density is
insensitive to changes in the packing in any bounded region. For
example, a finite number of balls can be removed from the
face-centered cubic packing without affecting its density.

Consequently, it is not possible to hope for any strong uniqueness
results for packings of optimal density.  The uniqueness
established by Lemma \ref{lemma:unique} is nearly as strong as can
be hoped for. It shows that certain local structures
(decomposition stars) attached to the face-centered cubic (fcc)
and hexagonal-close packings (hcp) are the only structures that
maximize a local density function.

\subsection{Hariot and Kepler}

\vskip -5mm \hspace{5mm}

 The modern mathematical study of close
packings can be traced to T. Hariot.  Hariot's work---unpublished,
unedited, and largely undated---shows a preoccupation with
packings of balls. He seems to have first taken an interest in
packings at the prompting of Sir Walter Raleigh. At the time,
Hariot was Raleigh's mathematical assistant,  and Raleigh gave him
the problem of determining formulas for the number of cannonballs
in regularly stacked piles. Shirley, Hariot's biographer, writes
that this study ``led him inevitably to the corpuscular or atomic
theory of matter originally deriving from Lucretius and
    Epicurus \cite[p.242]{Shi83}.''

Kepler became involved in packings of balls through his
correspondence with Hariot in the early years of the 17th century.
Kargon writes, in his history of atomism in England,

{ \narrower \font\ninerm=cmr9 \ninerm

    According to Hariot the universe is composed of atoms with void space interposed. The atoms themselves are eternal and
    continuous. Physical properties result from the magnitude, shape, and motion of these atoms, or corpuscles compounded from them$\ldots$.

    Probably the most interesting application of Hariot's atomic theory was in the field of optics. In a letter to Kepler on 2 December 1606
    Hariot outlined his views. Why, he asked, when a light ray falls upon the surface of a transparent medium, is it partially reflected and
    partially refracted? Since by the principle of uniformity, a single point cannot both reflect and transmit light, the answer must lie in the
    supposition that the ray is resisted by some points and not
    others.$\ldots$

    It was here that Hariot advised Kepler to abstract himself mathematically into an atom in order to enter `Nature's house'. In his reply of 2
    August 1607, Kepler declined to follow Harriot, ad atomos et vacua. Kepler preferred to think of the reflection-refraction problem in terms
    of the union of two opposing qualities---transparence and opacity. Hariot was surprised. ``If those assumptions and reasons satisfy you, I
    am amazed.'' \cite[p.26]{Kar66}

}

\smallskip
Despite Kepler's initial reluctance to adopt an atomic theory, he
was eventually swayed, and in 1611 he published an essay that
explores the consequences of a theory of matter composed of small
spherical particles.  Kepler's essay was the ``first recorded step
towards a mathematical theory of the genesis of inorganic or
organic form'' \cite[p.v]{Why66}.

Kepler's essay describes the face-centered cubic packing and
asserts that ``the packing will be the tightest possible, so that
in no other arrangement  could more pellets be stuffed into the
same container.''  This assertion has come to be known as the
Kepler conjecture.   This conjecture was verified with computer
assistance in 1998 \cite{part6}.

\subsection{Newton and Gregory}

\vskip -5mm \hspace{5mm}

The next episode in the history of this problem is a debate
between Isaac Newton and David Gregory.
  Newton and Gregory discussed the question of how many balls
of equal radius can be arranged to touch a given ball.  This is
the three-dimensional analogue of the simple fact that in two
dimensions six pennies, but no more, can be arranged to touch a
central penny.  This is the kissing-number problem in
$n$-dimensions. In three dimensions, Newton said that the maximum
was 12 balls, but Gregory claimed that 13 might be possible.  B.
L. van der Waerden and Sch\"utte in 1953 showed that Newton was
correct \cite{Sch53}.

The two-dimensional analogue of the Kepler conjecture is to show
that the honeycomb packing in two dimensions gives the highest
density.  This result was established in 1892 by Thue, with a
second proof appearing in 1910 (\cite{Thu92}, \cite{Thu10}).

In 1900, Hilbert made the Kepler conjecture part of his 18th
problem \cite{Hil01}.  The third part of that problem asks, `` How
can one arrange most densely in space an infinite number
    of equal solids of given form, e.g. spheres with given radii $\ldots$,
    that is, how can one so fit them together that the ratio of the
    filled to the unfilled space may be as great as possible?''

\subsection{The literature}

\vskip -5mm \hspace{5mm}

 Past progress toward the Kepler
conjecture can be arranged into four categories: (1) bounds on the
density, (2) descriptions of classes of packings for which the
bound of $\pi/\sqrt{18}$ is known, (3) convex bodies other than
balls for which the packing density can be determined precisely,
(4) strategies of proof.

Various upper bounds have been established on the density of
packings.  A list of such bounds appears in \cite{historical}.
Rogers's bound of $0.7797$ is particularly natural \cite{Rog58}.
It remained the best available bound for many years.

\subsection{Classes of packings}

\vskip -5mm \hspace{5mm}

 If the infinite dimensional space of all
packings is too unwieldy, we can ask if it is possible to
establish the bound $\pi/\sqrt{18}$ for packings with special
structures.

If we restrict the problem to packings of balls whose centers are
the points of a lattice, the
 packings are described by a finite number of parameters, and the
problem becomes much more accessible.  Lagrange proved that the
densest lattice packing in two dimensions is the familiar
honeycomb arrangement \cite{Lag73}. Gauss proved that the densest
lattice packing in three dimensions is the face-centered cubic
\cite{Gau31}.   The enormous list of references in \cite{CoSl98}
documents the many developments in lattice packings over the past
two centuries.

\subsection{Other convex bodies}

\vskip -5mm \hspace{5mm}

 If the optimal packings of balls are too
difficult to determine, we might ask whether the problem can be
solved for other convex bodies. To avoid trivialities, we restrict
our attention to convex bodies whose packing density is strictly
less than 1.

  The first convex body in Euclidean 3-space that does not tile
for which the packing density was explicitly determined is an
infinite cylinder \cite{Bez90}. Here A. Bezdek and W. Kuperberg
prove that the optimal density is obtained by arranging the
cylinders in parallel columns in the honeycomb arrangement.

In 1993, J. Pach exposed the humbling depth of our ignorance when
he issued the challenge to determine the packing density for some
bounded convex body that does not tile space \cite{MP93}.  (This
challenge was met  by A. Bezdek \cite{Bez94}.)

\subsection{Strategies of proof}

\vskip -5mm \hspace{5mm}

In 1953, L. Fejes T\'oth proposed a program to prove the Kepler
conjecture \cite{Fej53}. A single Voronoi cell cannot lead to a
bound better than the dodecahedral bound.  (The dodecahedral bound
is the ratio of the volume of a inscribed ball to the volume of
the containing dodecahedron.)  L. Fejes T\'oth considered weighted
averages of the volumes of collections of Voronoi cells.
 These weighted
averages involve up to 13 Voronoi cells.  He showed that if a
particular weighted average of volumes is greater than the volume
of the rhombic dodecahedron, then the Kepler conjecture follows.
The Kepler conjecture is an optimization problem in an infinite
number of variables.  L. Fejes T\'oth's weighted-average argument
was the first indication that it might be possible to reduce the
Kepler conjecture to a problem in a finite number of variables.
Needless to say, calculations involving the weighted averages of
the volumes of several Voronoi cells are complex.

L. Fejes T\'oth made another significant suggestion in
\cite{Fej64}. He was the first to suggest the use of computers in
the Kepler conjecture. After describing his program, he writes,

{ \narrower \font\ninerm=cmr9 \ninerm

Thus it seems that the problem can be reduced to the determination
of the minimum of a function of a finite number of variables,
providing a programme realizable in principle.  In view of the
intricacy of this function we are far from attempting to determine
the exact minimum.  But, mindful of the rapid development of our
computers, it is imaginable that the minimum may be approximated
with great exactitude.

}

\smallskip
A widely publicized attempt to prove the Kepler conjecture was
that of Wu-Yi Hsiang \cite{Hsi93}, \cite{Hsi02}.  Hsiang's
approach can be viewed as a continuation and extension of L. Fejes
T\'oth's program. Hsiang's work contains major gaps and errors
\cite{CoHMS94}.  A list of published materials relating to these
errors can be found in \cite{historical}.

\section*{2. Structure of the proof} \addsec
\label{sec:overview}

\vskip -5mm \hspace{5mm}

This section describes the structure of the proof of the Kepler
Conjecture.

\begin{theorem}
\label{theorem:kepler} (The Kepler Conjecture)  No packing of
congruent balls in Euclidean three space has density greater than
that of the face-centered cubic packing.
\end{theorem}

Here, we describe the top-level outline of the proof and give
references to the sources of the details of the proof
(\cite{formulation}, \cite{part1}, \cite{part2}, \cite{part3},
\cite{part4}, \cite{thesis}, \cite{part6}).

Consider a packing of congruent balls of unit radius in Euclidean
three space. The density of a packing does not decrease when balls
are added to the packing. Thus, to answer a question about the
greatest possible density we may add non-overlapping balls until
there is no room to add further balls. Such a packing will be said
to be {\it saturated}.

Let $\Lambda$ be the set of centers of the balls in a saturated
packing.  Our choice of radius for the balls implies that any two
points in $\Lambda$ have distance at least $2$ from each other. We
call the points of $\Lambda$ {\it vertices}.  Let $B(x,r)$ denote
the ball in Euclidean three space at center $x$ and radius $r$.
Let $\delta(x,r,\Lambda)$ be the finite density, defined by the
ratio of $A(x,r,\Lambda)$ to the volume of $B(x,r)$, where
$A(x,r,\Lambda)$ is defined as the volume of the intersection with
$B(x,r)$ of the union of all balls in the packing. Set
$\Lambda(x,r) = \Lambda \cap B(x,r)$.

The {\it Voronoi cell\/} $\Omega(v)$ around a vertex $v\in
\Lambda$ is the set of points closer to $v$ than to any other ball
center. The volume of each Voronoi cell in the face-centered cubic
packing is $\sqrt{32}$.  This is also the volume of each Voronoi
cell in the hexagonal-close packing.

Let $a:\Lambda\to\R$ be a function.  We say that $a$ is {\it
negligible\/} if there is a constant $C_1$ such that  for all
$x\in \R^3$ and $r\ge1$, we have
$$\sum_{v\in\Lambda(x,r)} a(v) \le C_1 r^2.$$
We say that the function $a$ is {\it fcc-compatible\/} if for all
$v\in\Lambda$ we have the inequality
$$\sqrt{32}\le \Vol(\Omega(v)) + a(v).$$

\begin{lemma}
\label{lemma:deltabound} If there exists a negligible
fcc-compatible function $a:\Lambda\to\R$ for a saturated packing
$\Lambda$, then there exists a constant $C$ such that for all
$x\in \R^3$ and $r\ge1$, we have
$$\delta(x,r,\Lambda)
\le \pi/\sqrt{18} + C/r.$$
\end{lemma}

\noindent{\bf Proof.}  The numerator $A(x,r,\Lambda)$ of
$\delta(x,r,\Lambda)$ is at most the product of the volume of a
ball $4\pi/3$ with the number $|\Lambda(x,r+1)|$ of balls
intersecting $B(x,r)$.  Hence
    \begin{equation}
    A(x,r,\Lambda) \le |\Lambda(x,r+1)| 4\pi/3.
    \label{eqn:Abound}
    \end{equation}

In a saturated packing each Voronoi cell is contained in a ball of
radius $2$ centered at the {\it center} of the cell.  The volume
of the ball $B(x,r+3)$ is at least the combined volume of Voronoi
cells lying entirely in the ball. This observation, combined with
fcc-compatibility and negligibility, gives
    \begin{equation}
    \begin{split}
    \sqrt{32}|\Lambda(x,r+1)|
    &\le \sum_{v\in\Lambda(x,r+1)} (a(v) +
    \Vol(\Omega(v))) \\
    &\le C_1 (r+1)^2 + \Vol\,B(x,r+3) \\
    &\le C_1 (r+1)^2 + (1+3/r)^3 \Vol\,B(x,r).
    \label{eqn:Bbound}
    \end{split}
    \end{equation}
Divide through by $\Vol\,B(x,r)$ and eliminate $|\Lambda(x,r+1)|$
between Inequality (\ref{eqn:Abound}) and Inequality
(\ref{eqn:Bbound}) to get
    $$\delta(x,r,\Lambda)
        \le \frac{\pi}{\sqrt{18}} (1+3/r)^3 + C_1 \frac{(r+1)^2}{r^3\sqrt{32}}.
    $$
The result follows for an appropriately chosen constant $C$.

\begin{remark} We take the precise meaning of the Kepler Conjecture to
be a bound on the essential supremum of the function $\delta(x,r)$
as $r$ tends to infinity. Lemma \ref{lemma:deltabound} implies
that the essential supremum of $\delta(x,r,\Lambda)$ is bounded
above by $\pi/\sqrt{18}$, provided a negligible fcc-compatible
function can be found.  The strategy will be to define a
negligible function, and then to solve an optimization problem in
finitely many variables to establish that it is fcc-compatible.
\end{remark}

The article \cite{formulation} defines a compact topological space
$X$ and a continuous function $\sigma$ on that space.

The topological space $X$ is directly related to packings. If
$\Lambda$ is a saturated packing, then there is a geometric object
$D(v,\Lambda)$ constructed around each vertex $v\in\Lambda$.
$D(v,\Lambda)$ depends on $\Lambda$ only through the vertices in
$\Lambda$ at distance at most $4$ from $v$.  The objects
$D(v,\Lambda)$ are called {\it decomposition stars}, and the space
of all decomposition stars is precisely $X$.

Let $\dtet$ be the packing density of a regular tetrahedron.  That
is, let $S$ be a regular tetrahedron of edge length $2$.  Let $B$
the part of $S$ that lies within distance $1$ of some vertex. Then
$\dtet$ is the ratio of the volume of $B$ to the volume of $S$. We
have $\dtet = \sqrt{8} \arctan(\sqrt{2}/5)$.

Let $\doct$ be the packing density of a regular octahedron of edge
length $2$, again constructed as the ratio of the volume of points
within distance $1$ of a vertex to the volume of the octahedron.
We have $\doct \approx 0.72$.

Let $\pt = -\pi/3 + \sqrt{2}\dtet\approx 0.05537$.

The following conjecture is made in \cite{formulation}:

\begin{conjecture}
\label{conjecture:sigma}
 The maximum of $\sigma$ on $X$ is the constant
$8\,\pt\approx 0.442989$.
\end{conjecture}

\begin{lemma}
\label{lemma:exista} An affirmative answer to Conjecture
\ref{conjecture:sigma} implies the existence of a negligible
fcc-compatible function for every saturated packing $\Lambda$.
\end{lemma}

\noindent{\bf Proof.}  For any saturated packing $\Lambda$ define
a function $a:\Lambda\to\R$ by
$$-\sigma(D(v,\Lambda))/(4\doct) + 4\pi/(3\doct) = \Vol(\Omega(v)) + a(v).$$
Negligibility follows from \cite[Prop. 3.14 (proof)]{formulation}.
The upper bound of $8\,\pt$ gives a lower bound $$-8\,\pt/(4\doct)
+ 4\pi/(3\doct) \le \Vol(\Omega(v)) + a(v).$$ The constant on the
left-hand side of this inequality equals $\sqrt{32}$, and this
establishes fcc-compatibility.

\begin{theorem}
\label{theorem:sigma} Conjecture \ref{conjecture:sigma} is true.
That is,
  the maximum of the function $\sigma$ on the
topological space $X$ of all decomposition stars is $8\,\pt$.
\end{theorem}

Theorem \ref{theorem:sigma}, Lemma \ref{lemma:exista}, and Lemma \ref{lemma:deltabound} combine to give a proof of
the Kepler Conjecture \ref{theorem:kepler}.

Let $t_0=1.255$ ($2t_0 = 2.51$).  This is a parameter that is used
for truncation throughout the series of articles on the Kepler
Conjecture.

Let $U(v,\Lambda)$ be the set of vertices in $\Lambda$ at distance
at most $2t_0$ from $v$.  From a decomposition star $D(v,\Lambda)$
it is possible to recover $U(v,\Lambda)$ (at least up to Euclidean
translation: $U\mapsto U+y$, for $y\in\R^3$).  We can completely
characterize the decomposition stars at which the maximum of
$\sigma$ is attained.

\begin{theorem}
\label{theorem:sharp} Let $D$ be a decomposition star at which the
maximum $8\,\pt$ is attained.  Then the set $U(D)$ of vectors at
distance at most $2t_0$ from the center has cardinality $12$. Up
to Euclidean motion, $U(D)$ is the kissing arrangement of the $12$
balls around a central ball in the face-centered cubic packing or
hexagonal-close packing.
\end{theorem}

\subsection{Outline of proofs}

\vskip -5mm \hspace{5mm}

To prove Theorems \ref{theorem:sigma} and \ref{theorem:sharp}, we
wish to show that there is no counterexample.  That is, we wish to
show that there is no decomposition star $D$ with value $\sigma(D)
> 8\,\pt$.  We reason by contradiction, assuming the existence of
such a decomposition star.  With this in mind, we call $D$ a {\it
contravening decomposition star}, if
    $$\sigma(D)\ge 8\,\pt.$$
In much of what follows we will assume that every decomposition
star under discussion is a contravening one.  Thus, when we say
that no decomposition stars exist with a given property, it should
be interpreted as saying that no such contravening decomposition
stars exist.

To each contravening decomposition star, we associate a
(combinatorial) plane graph.  A restrictive list of properties of
plane graphs is described in \cite[Section 2.3]{part6}.  Any plane
graph satisfying these properties is said to be {\it tame}. All
tame plane graphs have been classified. (There are several
thousand, up to isomorphism.) Theorem \cite[Th 2.1]{part6} asserts
that the plane graph attached to each contravening decomposition
star is tame. By the classification of such graphs, this reduces
the proof of the Kepler Conjecture to the analysis of the
decomposition stars attached to the finite explicit list of tame
plane graphs.

A few of the tame plane graphs are of particular interest. Every
decomposition star attached to the face-centered cubic packing
gives the same plane graph (up to isomorphism).  Call it
$G_{fcc}$.  Likewise, every decomposition star attached to the
hexagonal-close packing gives the same plane graph $G_{hcp}$. Let
$X_{crit}$ be the set of decomposition stars $D$ such that the set
$U(D)$ of vertices is the kissing arrangement of the $12$ balls
around a central ball in the face-centered cubic or
hexagonal-close packing.  There are only finitely many orbits of
$X_{crit}$ under the group of Euclidean motions.
\begin{figure}[htb]
  \centering
  \includegraphics{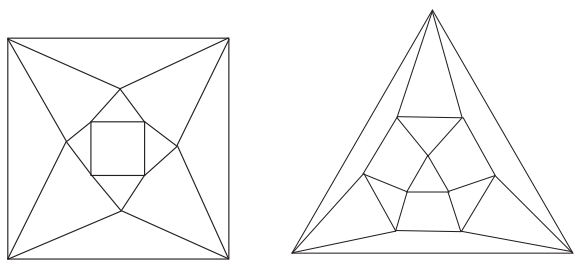}
  \caption{The plane graphs $G_{fcc}$ and $G_{hcp}$}
  \label{fig:gfcchcp}
\end{figure}

In \cite[Lemma 3.13]{formulation}, the necessary local analysis is
carried out to prove the following local optimality.

\begin{lemma} \label{lemma:unique}
A decomposition star whose plane graph is $G_{fcc}$ or $G_{hcp}$
has score at most $8\,\pt$, with equality precisely when the
decomposition star belongs to $X_{crit}$. \end{lemma}

In light of this result, we prove \ref{theorem:sigma} and
\ref{theorem:sharp} by proving that any decomposition star whose
graph is tame and not equal to $G_{fcc}$ or $G_{hcp}$ is not
contravening

\begin{figure}[htb]
  \centering
  \includegraphics{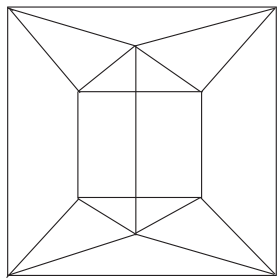}
  \caption{The plane graph $G_{pent}$}
  \label{fig:pentahedral}
\end{figure}

There is one more tame plane graph that is particularly
troublesome.  It is the graph $G_{pent}$ obtained from the
pictured configuration of twelve balls tangent to a given central
ball (Figure \ref{fig:pentahedral}). (Place a ball at the north
pole, another at the south pole, and then form two pentagonal
rings of five balls.) This case requires individualized attention.
S. Ferguson proves in \cite{thesis} that if $D$ is any
decomposition star with this graph, then $\sigma(D)< 8\,\pt$.

To eliminate the remaining cases, more-or-less generic arguments
can be used.  A linear program is attached to each tame graph $G$.
The linear program can be viewed as a linear relaxation of the
nonlinear optimization problem of maximizing $\sigma$ over all
decomposition stars with a given tame graph $G$. Because it is
obtained by relaxing the constraints on the nonlinear problem, the
maximum of the linear problem is an upper bound on the maximum of
the original nonlinear problem.  Whenever the linear programming
maximum is less than $8\,\pt$, it can be concluded that there is
no contravening decomposition star with the given tame graph $G$.
This linear programming approach eliminates most tame graphs.

When a single linear program fails to give the desired bound, it
is broken into a series of linear programming bounds, by branch
and bound techniques.  For every tame plane graph $G$ other than
$G_{hcp}$, $G_{fcc}$, and $G_{pent}$, we produce a series of
linear programs that establish that there is no contravening
decomposition star with graph $G$.  When every face of the plane
graph is a triangle or quadrilateral, this is accomplished in
\cite{part3}.  The general case is completed in the final sections
of \cite{part6}.

\label{lastpage}


\begin{thebibliography}{aa}


\bibitem{Bez90} A. Bezdek and W. Kuperberg, Maximum density space
packing with
    congruent circular cylinders of infinite length,
    {\it Mathematica} 37 (1990), 74--80.

\bibitem{Bez94} A. Bezdek, A remark on the packing density in the
3-space
    in {\it Intuitive Geometry}, ed. K. B\"or\"oczky and G. Fejes
    T\'oth, {\it Colloquia Math. Soc. J\'anos Bolyai} 63, North-Holland
    (1994), 17--22.






\bibitem{CoHMS94} J. H. Conway, T. C. Hales, D. J. Muder, and N. J. A.
Sloane,
    On the Kepler conjecture, {\it Math. Intelligencer} 16,
    no. 2 (1994), 5.

\bibitem{CoSl98} J. H. Conway and N. J. A. Sloane, Sphere packings,
lattices
    and groups,  third edition, Springer-Verlag, New York,
    1998.


\bibitem{Fej53} L. Fejes T\'oth, {\it Lagerungen in der Ebene auf der
Kugel und im Raum}, Springer, Berlin, first edition, 1953.

\bibitem{Fej64} L. Fejes T\'oth, Regular figures, Pergamon Press,
    Oxford London New York, 1964.

\bibitem{thesis} S. Ferguson, Sphere Packings V, thesis, University of
Michigan,
    1997.

\bibitem{formulation} S. Ferguson, T. Hales, A Formulation of the
Kepler Conjecture, preprint math.MG/9811072.


\bibitem{Gau31} C. F. Gauss, Untersuchungen \"uber die Eigenscahften der
positiven tern\"aren quadratischen Formen von Ludwig August Seber,
    {\it G\"ottingische gelehrte Anzeigen}, 1831 Juli 9,
also in {\it J. reine angew. Math.} 20 (1840), 312--320.





\bibitem{historical} T. C. Hales, An Overview of the Kepler
Conjecture, preprint math.MG/9811071.

\bibitem{part1} T. C. Hales, Sphere Packings I,
    Discrete and Computational Geometry, 17 (1997), 1-51.

\bibitem{part2} T. C. Hales, Sphere Packings II,
    Discrete and Computational Geometry, 18 (1997), 135-149.

\bibitem{part3} T. C. Hales, Sphere Packings III, preprint math.MG/9811075.

\bibitem{part4} T. C. Hales, Sphere Packings IV, preprint math.MG/9811076.

\bibitem{part6} T. C. Hales, The Kepler Conjecture, preprint math.MG/9811078.







\bibitem{Hil01} D. Hilbert, Mathematische Probleme, {\it Archiv Math.
Physik} 1 (1901),
    44--63, also in {\it Proc. Sym. Pure Math.} 28 (1976), 1--34.



\bibitem{Hsi93} W.-Y. Hsiang, On the sphere packing problem and the
proof of Kepler's conjecture, Internat. J. Math 93 (1993),
739-831.

\bibitem{Hsi02} W.-Y. Hsiang, Least Action Principle of Crystal
Formation of Dense Packing Type and the Proof of Kepler's
Conjecture, World Scientific, 2002.





\bibitem{Kar66} R. Kargon, Atomism in England from Hariot to Newton,
    Oxford, 1966.

\bibitem{Kep66} J. Kepler, The Six-cornered snowflake, Oxford Clarendon
Press,
    Oxford, 1966,  forward by L. L. Whyte.

\bibitem{Lag73} J. L. Lagrange,  Recherches d'arithm\'etique, {\it Nov.
Mem.
    Acad. Roy. Sc. Bell Lettres Berlin} 1773, in {\it \OE uvres}, vol. 3,
    693--758.



\bibitem{MP93} W. Moser, J. Pach, Research problems in discrete
geometry,
    DIMACS Technical Report, 93032, 1993.





\bibitem{Rog58} C. A. Rogers, The packing of equal spheres, {\it Proc.
London Math.
    Soc.} (3) 8 (1958), 609--620.

\bibitem{Sch53} K. Sch\"utte and B.L. van der Waerden, Das Problem der
dreizehn Kugeln, {\it Math. Annalen} 125, (1953), 325--334.

\bibitem{Shi83} J. W. Shirley, {\it Thomas Harriot: a biography},
Oxford, 1983.

\bibitem{Thu92} A. Thue, Om nogle geometrisk taltheoretiske Theoremer,
    {\it Forandlingerneved de Skandinaviske Naturforskeres} 14 (1892), 352--353.

\bibitem{Thu10} A. Thue, \"Uber die dichteste Zusammenstellung von
    kongruenten Kreisen in der Ebene, {\it Christinia Vid. Selsk. Skr.} 1
    (1910), 1--9.

\bibitem{Why66} L. L. Whyte, forward to \cite{Kep66}.

\end{thebibliography}
\end{document}